\documentclass[12pt]{article}

\usepackage{latexsym}
\usepackage{amsfonts}
\usepackage{amssymb}
\usepackage{amsmath}
\usepackage{mathrsfs}
\usepackage[all]{xy}
\usepackage{graphicx}
\usepackage{stmaryrd}

\newtheorem{Theory}{Theory}[subsection] 
\newtheorem{theorem}[Theory]{Theorem}
\newtheorem{lemma}[Theory]{Lemma}

\newtheorem{corollary}[Theory]{Corollary}

\newtheorem{remark}[Theory]{Remark} 
\newtheorem{claim}[Theory]{Claim}

\newcommand{\Z}{\mathbb{Z}}

\newcommand{\N}{\mathbb{N}}

\newcommand{\ploi}{\mbox{PL}_0(\mbox{{\bf{I}}})}

\newcommand{\be}{\begin{enumerate}}
\newcommand{\ee}{\end{enumerate}}
\newcommand{\bt}{\begin{theorem}}
\newcommand{\et}{\end{theorem}}
\newcommand{\bl}{\begin{lemma}}
\newcommand{\el}{\end{lemma}}
\newcommand{\bc}{\begin{corollary}}
\newcommand{\ec}{\end{corollary}}
\newcommand{\br}{\begin{remark}}
\newcommand{\er}{\end{remark}}
\newcommand{\bcl}{\begin{claim}}
\newcommand{\ecl}{\end{claim}}

\newcommand{\supp}{\operatorname{Supp}}

\begin{document}

\begin{center}
{\bf{Subgroups of R. Thompson's Group $F$ that are Isomorphic to $F$.}}

\vspace{8pt}

Bronlyn Wassink

\end{center}

\vspace{10pt} {\flushleft ABSTACT: }

This paper studies when a pair of elements in $F$ are the images of the standard
generators of F under a self monomorphism.

\vspace{15pt}
\section{Introduction}

Richard Thompson's Groups $F \subseteq T \subseteq V$ were defined by Richard
Thompson in 1965.  Thompson proved that both $T$ and $V$ are finitely
presented infinite simple groups.  The group $F$ is also finitely presented,
but it is not simple.  A good introduction to these groups is in \cite{CFP}.

In this paper we find all subgroups of $F$ isomorphic to $F$ by finding each
pair or elements of $F$ that is the image, under a monomorphism of $F$ into
itself, of a standard pair of generaotrs of $F$.  This paper was motivated by the questions about the distortion of subgroups of
$F$.  It has been shown by Jos\'{e} Burillo in \cite{BQI} that for any natural
number $n$, there is a subgroup of $F$ that is isomorphic to $F^n \times \Z$
that can be quasi-isometrically embedded into $F$ using the inclusion
function.  It is unknown if every subgroup of $F$ that is isomorphic to $F$ or
$F \times \Z$ has this property.  In order get started on this question, we
must first find every subgroup of $F$ that is isomorphic to $F$.

There are several ways to define Thompson's Groups.  The definition that this
paper will focus on defines $F$ as a subgroup of $\ploi$, which is the group
of all piecewise-linear, orientation preserving homeomorphisms of the unit
interval, {\textbf{I}}, that admit only finitely many points of
non-differentiability.  The group $F$ is exactly the elements of $\ploi$
whose points of non-differentiability all of occur at a dyadic
rationals and have all slopes integral powers of $2$.  The dyadic
rationals are defined to be the rational numbers of the form $\frac{m}{2^n}$,
where $m, n \in \Z$ and $n>0$.

In the first section of this paper, many standard results and terminology
will be introduced.

The second section of this paper gives our main result.   
The presentation $\langle x_0, x_1 | [x_0x_1^{-1}, x_1^{x_0}] = [x_0x_1^{-1},
x_1^{x_0}] =1  \rangle$ is called the ``standard'' finite presentation for
$F$.  We find all pairs $(f,g)$ of elements in $\ploi$ for which there is an isomorphic embedding of $F$ into $\ploi$ that carries $x_0$ to $f$ and $x_1$ to $g$.

The third section of this paper gives a partial converse to the main theorem in \cite{BrinU}.  This section also provides a counterexample to the full converse of the main result in \cite{BrinU}.
 
In addition to the application discussed in Section 3, the results in Section 2 of
this paper are used in the complete classification of
the isomorphism classes of the finite index subgroups of $F$.  This 
classification appears in \cite{BronlynCollin}, which is a joint paper of the
author and Collin Bleak.  Independently from our work, Burillo, Cleary and
R\"{o}ver in \cite{BCRCommensurator} also characterized the finite index
subgroups of $F$ that are isomorphic to $F$ as a result of their
investigations into the commensurator of $F$.  Their techniques are very
different from those in \cite{BronlynCollin}

\subsection{Definitions and Notation}

Richard Thompson's Group $F$ can also described by the following
presentations.

\begin{equation} \label{ei}
 F \cong \langle x_0, x_1, x_2, ... \ | \ x_j^{x_i} = x_{j+1} \ \mbox{for} \ i<j \rangle
\end{equation}

\begin{equation} \label{ef}
F  \cong \langle x_0, x_1 \ | \ [x_0x_1^{-1}, x_1^{x_0}]=[x_0x_1^{-1},
 x_1^{x_0^2}]=1   \rangle
\end{equation}

\noindent
where $a^{b} = b^{-1}ab$ and $[a,b] = aba^{-1}b^{-1}$.

Define the {\em{standard finite presentation of $F$}} to be the presentation
(\ref{ef}) above.

We will make extensive use of the representation of $F$ in the group $\ploi$,
whose definition is given in the introduction.

Composition and evaluation of functions in $F$ will
be in word order.  That is, if $f,g \in F$ and $t \in
[0,1]$,  then $tf = f(t)$, $fg = g \circ f$, and $f^g = g^{-1}fg = g
\circ f \circ g^{-1} $.

The functions $x_0$ and $x_1$ are defined below. The standard
representation of $F$ as PL homeomorphisms of {\bf{I}} carries $x_0$ and $x_1$
of (\ref{ef}) to the functions below with the same names.

\[ ax_0 = \left\{
\begin{array}{cl}
   2a   &   0 \leq a \leq 1/4      \\
   a + 1/4   &  1/4 \leq a \leq 1/2       \\
  a/2 + 1/2    &  1/2 \leq a \leq 1       \\
\end{array} \right. \ \ \ \ \   ax_1 = \left\{
\begin{array}{cl}
 a & 0 \leq a \leq 1/2 \\
   2a - 1/2   &   1/2 \leq a \leq 5/8      \\
   a + 1/8   &  5/8 \leq a \leq 3/4       \\
  a/2 + 1/2    &  3/4 \leq a \leq 1       \\
\end{array} \right. \ \ \ \ \    \]

One can check, using the conventions above, that the functions $x_0$ and
$x_1$ satisfy the relations in (\ref{ef}).  Presentations (\ref{ei}) and (\ref{ef}) are equivalent (see \cite{CFP} Theorem 3.4).  The fact that the functions $x_0$ and
$x_1$ generate all of the claimed functions in $F$ (as a group of
homeomorphisms) is Corollary 2.6 in \cite{CFP}.
(Note that our functions $x_0$ and $x_1$ are the inverses of the
homeomorphisms in \cite{CFP}.)

Given a homeomorphism $f:[0,1]\to [0,1]$, we will denote by $\supp(f)$
the {\em{support of $f$}}, defined as
\[
\supp(f)=\left\{x\in [0,1]|xf\neq f\right\}.
\]

Given any $f \in \ploi$, $\supp(f)$ will be a finite union of disjoint open intervals.  Each of these open intervals will be called an {\em{orbital}} of $f$.

As $F \leq \ploi$, many lemmas and theorems will be done using more general functions from $\ploi$ instead of the more restrictive functions in $F$.
 
We will say that two functions $f_0$ and $f_1$ in $\ploi$ generate a
{\em{standard isomorphic copy of $F$}} if the subgroup of $\ploi$ that is
generated by $f_0$ and $f_1$ is isomorphic to $F$ by an isomorphism sending
$x_0$ to $f_0$ and $x_1$ to $f_1$.  For $f_0$ and $f_1$ to generate a standard
isomorphic copy of $F$, it is necessary that $[f_0f_1^{-1},
f_1^{f_0}]=[f_0f_1^{-1}, f_1^{f_0^2}]=1 $.

\subsection{Previous Results}

The first remark below, though trivial, is useful is several proofs.  It follows
because the relations in the remark are relations in (\ref{ef}).

\br \label{infinite} If $f_0$ and $f_1$ are functions in $\ploi$ that satisfy both $[f_0f_1^{-1} , f_1^{f_0} ]=1$ and $[ f_0f_1^{-1} , f_1^{f_0^2} ]=1$, then there is a group homomorphism from $F$ to $\langle f_0, f_1 \rangle$ sending $x_0$ to $f_0$ and $x_1$ to $f_1$.   In particular, $f_0$ and $f_1$ will satisfy all of the relations that $x_0$ and $x_1$ they will satisfy in $F$, such as  $[f_1f_2^{-1} , f_2^{f_1} ]=1$, where $f_2 = f_1^{f_0}$.
 \er

\bl\label{disjointSupport} 

If $f,g \in \ploi$ where the support of $f$ is disjoint from the support of $g$, then $f$ and $g$ commute.

 \el

\bl \label{unionSupports}

Let $g_1, g_2, ..., g_n \in \ploi$.  Let $H$ be the subgroup of $\ploi$ that is generated
by $g_1, g_2, ..., g_n$ and define
\[
\supp(H) = \left\{x\in [0,1]|xh\ne x \textrm{ for some }h\in
H\right\}.
\]
Then, $\displaystyle \supp(H) = \bigcup_{i = 1}^{n} \supp(g_i)$.

\el

\vspace{8pt}

Whenever $H$ is a subgroup of $\ploi$, the support of $H$ will consist of disjoint open intervals.  Call each of those intervals an {\em{orbital of $H$}}.

The following Lemma is a restatement of Remark 2.1, Lemma 2.2 and Lemma 2.5 in  \cite{ansc}.

\bl\label{support}

\be
\item Any element $h \in \ploi$ has only finitely many orbitals.
\item Let $h$ be a function in $\ploi$.   If $A$ is an orbital for $h$, then either $xh>x$ or $xh<x$ for all points $x$ in $A$.
\item If $h \in \ploi$ and $A = (a,b)$ is an orbital of $h$, then given any $\varepsilon > 0$ and $x$ in $A$, there is an integer $n$ so that $|xh^{-n} - a| < \varepsilon$ and $|xh^{n} - b| < \varepsilon$.
\item Let $g,h \in \ploi$.  $(a,b)$ is an orbital of $g$ if and only if $(ah,bh)$ is an orbital of $g^h$.
\item If $(a,b)$ is an orbital of $f$, $c \in (a,b)$, and $m \in \Z$, then $cf^m \in (a,b)$.
\ee
\el

Lemma \ref{support} allows us to introduce new terminology.  From Lemma \ref{support}.(1), since there are only finitely many orbitals, we can introduce an order on them.  If $A$ and $B$ are orbitals of $g$, then we say that $A < B$ when $a<b$ for any $a \in A$ and $b \in B$.  With this ordering, we can refer to the {\em{first orbital}}, the {\em{second orbital}}, and so on.  Many times we will be interested in the first orbital and the last orbital.

Lemma \ref{support}.(2) allows us to introduce a sign on the orbitals.  If an
orbital $A$ of $g$ has the property that $ah > a$ for all $a$ in $A$, then we
call $A$ an {\em{up-bump}} of $g$.  Similarly, if an orbital $B$ of $g$ has
the property that $bh < b$ for all $b$ in $B$, then we call $B$ a
{\em{down-bump}} of $g$.  Notice that Lemma \ref{support}.(2) guarentees that every orbital is either an up-bump or a down-bump.

We call the orbitals $(a,b)$ of $g$ and $(ah, bh)$ of $g^h$ {\em{corresponding orbitals}}.

The following lemma is a part of Lemma 2.5 in \cite{ansc}.

\bl \label{conjslopes} If $f, g \in
\ploi$ and $(a,b)$ is an orbital of $f$, then the derivative from the right of $f$ at $a$ equals the derivative from the right of $f^g$ at $ag$ and the derivative from the left of $f$ at $b$ equals the derivative from the left of $f^g$ at $bg$. 
\el

In other words, $f^{g}$ has the same {\em{leading}} and {\em{trailing slopes}} on its orbitals as $f$ has on each of its corresponding orbitals. 

\br \label{genslopes}
Given any functions $f_0$ and $f_1$ in $\ploi$ and $k \geq 1$, we can define
$f_{k+1} = f_1^{f_0^k}$.  So then all of the $f_k$ for $k \geq 1$ have the
same leading slopes on each of their corresponding orbitals, and all $f_k$
have the same trailing slopes on their corresponding orbitals.

\er

The following is Theorem 4.3 from \cite{CFP}.

\bl \label{npnaq}

$F$ has no proper non-abelian quotients. 

\el

In particular, if we can find an $f$ and $g$ where
\[
[fg^{-1},g^{f}]=[fg^{-1}, g^{f^2}]=1
\]
 and $f$ and $g$ do not commute, then $f$ and $g$ generate a group
that is isomorphic to $F$. 

The following lemma is a combination of Theorems 4.1 and 4.5 in \cite{CFP}.

\bl \label{FPrime}

The group $F' = [F,F]$, the commutator subgroup of $F$, is simple.
Furthermore, $\supp(F') = (0,1)$ and $F'$ consists of all of the functions $f \in F$ such that
both $f'(0) = 1$ and $f'(1) = 1$, where $f'(0)$ is the derivative of $f$ at zero from the right and  $f'(1)$ is the derivative of $f$ at one from the left.

\el

\vspace{8pt}

The following is a result by Brin and Squier in Theorem 4.18 in \cite{CommonRoot}.

\bt \label{commonroot} If $A$ is an orbital of both $f$ and $g$ in $\ploi$ and $[f,g]|_{A} = 1$, then there is a function $h$ and integers $m$ and $n$ so that $h^m|_A = f|_A$ and $h^n|_A = g|_A$.

\et

\section{Standard Isomorphic Copies of $F$.}

In this section, we characterize those pairs $(f,g)$ of elements in $\ploi$ that generate a standard isomorphic copy of $F$.

This section will be organized as follows.  Section 2.1 gives information
about commuting functions.   Then, assuming that a pair $(f_0, f_1)$ of functions in $\ploi$ generate a standard isomorphic copy of $F$, Sections 2.2 and 2.3 accumulate a series of increasingly restrictive observations extracted from the fact that $f_0$ and $f_1$ must satisfy the relations in (2).  These give a set of conditions that are necessary for $(f_0, f_1)$ to generate a standard isomorphic copy of $F$, but they will not quite be sufficient.  They are sufficient however to imply that $f_0$ and $f_1$ satisfy the second relation in (2).  Sections 2.4 and 2.5 then find the remaining conditions needed to characterize when $(f_0, f_1)$ generate a standard isomorphic copy of $F$.

In the introduction we remark that for $f_0$ and $f_1$ in $\ploi$ to generate
a standard isomorphic copy of $F$ it is necessary that $[f_0f_1^{-1},
f_1^{f_0}] = [f_0f_1^{-1}, f_1^{f_0^2}] = 1$ hold.  By Lemma \ref{npnaq}, we
see that this is a sufficient condition when $f_0$ and $f_1$ do not commute.

If $f_0$ and $f_1$ in $\ploi$ generate a standard isomorphic copy of $F$, we define $f_{k+1} = f_1^{f_0^k}$ for all $k \geq 0$.

\subsection{Commuting Elements}

\bl \label{commutedisj}
Let $f,g \in \ploi$ where $[f,g] = 1$.  If $A$ is an
orbital of $f$, either $A$ is disjoint from the support of $g$ or $A$ is also
an orbital of $g$.

\el

{\em{Proof:}}  Suppose $fg=gf$.  Let $A = (a,c)$ be an orbital of $f$ that is
not disjoint from the support of $g$.  Let $B$ be an orbital of $g$ where $A
\not= B$ and $A \cap B \not= \emptyset$.   By these assumptions, either $f$
has a fixed point in $B$ or $g$ has a fixed point in $A$.  Without loss of
generality, we may assume that $g$ has a fixed point in $A$ and that $A$ is an
up-bump of $f$. 

Suppose that $y \in (a,c)$ is a fixed point of $g$. For all $k \in \Z$ we have $yf^kg = ygf^k = yf^k$.  So $yf^k$ is in the interval $(a,c)$ and is not in an
orbital of $g$.  Since $g$ can only have finitely many orbitals and
$\displaystyle \lim_{k \to -\infty} bf^k = a$, then there must be some
interval $(a,t)$ that is disjoint from the support of $g$.  Because there are
points in $(a,c)$ that are in the support of $g$, then there will be an
orbital $(b,d)$ of $g$ where $b \in (a,c)$ and $(a,b)$ is disjoint from the
support of $g$.  Let $x \in (b, fb) \cap (b,d)$.  Let $z = xf^{-1}$.  So $z$
is a fixed point of $g$ and $zf$ is not a fixed point of $g$.  Then $(zg)f =
zf \not= (zf)g$, a contradiction.   $\square$

\bl \label{eqcommute}  Let $f_0,f_1 \in \ploi$ \ such that $[f_1^{f_0},f_0f_1^{-1}] = 1$ and \ $[f_0f_1^{-1},f_1^{f_0^2}] = 1$.  If $A$ is an orbital of both $f_0$ and $f_1$, then $f_0$ and $f_1$ commute on $A$.

\el

{\em{Proof:}}  Suppose $A = (a,c)$ is an orbital of both $f_0$ and $f_1$.
Since $f_{k+1} = f_1^{f_0^k}$, $ (af_0^k, cf_0^k) = (a,c) = A$ is also an
orbital for all $f_{k+1}$, $k \geq 1$.  Note that $[f_1f_2^{-1},f_1^{-1}f_2f_1] =
1$.  By Remark \ref{genslopes}, the leading slopes on $A$ of $f_1$ and $f_2$
are equal and the trailing slopes on $A$ of $f_1$ and $f_2$ are equal. So
there is an $a' \not= a$ with $a \in A$ and a $c' \not= c$ with $c \in A$ such that
$f_1|_{[a,a']} = f_2|_{[a,a']}$ and  $f_1|_{[c',c]} = f_2|_{[c',c]}$.  Then
$f_1f_2^{-1}$ is the identity on $[a,a'] \cup [c',c]$.  So then any orbital of
$f_1f_2^{-1}$ is either disjoint from $[a,c]$ or is contained in $[a',c']$. 
Since $ f_3$ and $f_1f_2^{-1}$ commute, by Lemma \ref{commutedisj} the
orbitals of these functions must either be disjoint or equal. 
If $P \subseteq [a',c']$ is an orbital $f_1f_2^{-1}$, then $P$ and $(a,c)$ can
neither be disjoint nor equal.  So then $P = \emptyset \Rightarrow
(f_1f_2^{-1})|_{A} = 1 \Rightarrow f_1|_A = f_2|_A$.  Thus $
(f_0^{-1}f_1f_0)|_{A}=  f_2|_{A} = f_1|_{A} \Rightarrow (f_1f_0)|_{A} =
(f_0f_1)|_{A}$.   $\square$

\subsection{Orbital Arrangements}

\bl \label{nests}
Let $f_0,f_1 \in \ploi$ \ such that $[f_1^{f_0},f_0f_1^{-1}] = 1$ and \ $[f_0f_1^{-1},f_1^{f_0^2}] = 1$.    Let $A_1, A_2, ..., A_N$ be the orbitals of $f_0$.  Let $B_1, B_2, ..., B_M$ be the orbitals of $f_1$.   Then for any $A_j$ and $B_i$ where $A_j \cap B_i \not= \emptyset$, either $A_j \subseteq B_i$ or $B_i \subseteq A_j$.

\el

{\em{Proof:}}  Suppose there is an $A_j$ and a $B_i$ such that $A_j \cap B_i \not= \emptyset$ but $A_j \not\subseteq B_i$ and $B_i \not\subseteq A_j$. 
Let $A_j = (a,c)$ and $B_i = (b,d)$.  Either $a<b<c<d$ or $b<a<d<c$. 
 Assume that $a<b<c<d$. 
One of the orbitals of $f_1^{f_0}$ is $(bf_0, df_0)$. 
Since $b<c<d$, then $bf_0 < cf_0 < df_0$. 
Since $c$ is a fixed point of $f_0$, the previous inequality simplifies to $bf_0 < c < df_0$.  So $c$ is in an orbital of $f_1^{f_0}$.

Since $c$ is not a fixed point of $f_1$ and is a fixed point of $f_0$, \  $cf_0f_1^{-1} = cf_1^{-1} \not= c$.  So $c$ is not a fixed point of $f_0f_1^{-1}$.  So $c$ is in an orbital of $f_0f_1^{-1}$.

Since the orbital $(bf_0,df_0)$ of $f_1^{f_0}$ is not disjoint from the set of all orbitals of $f_0f_1^{-1}$ ($c$ is in their intersection), then, by Lemma \ref{commutedisj}, an orbital for $f_0f_1^{-1}$ must also equal $(bf_0,df_0)$.

The point $c$ is also in the interval $(bf_0^2,df_0^2)$, which is an orbital for $f_1^{f_0^2}$. 
Since the orbitals of $f_0f_1^{-1}$ and $f_1^{f_0^2}$ are not disjoint (again, both contain the point $c$), then, by Lemma \ref{commutedisj}, an orbital for $f_0f_1^{-1}$ must equal $(bf_0^2, df_0^2)$.

Since $b$ is not a fixed point for $f_0$, then either $bf_0<bf_0^2<c$ or $bf_0^2 < bf_0 <c$.  So $(bf_0^2,df_0^2)$ and $(bf_0,df_0)$ are non-equal, non-disjoint intervals.  Thus they can not both be orbitals for the same function, $f_0f_1^{-1}$, which is a contradiction.

Similarly, if $b<a<d<c$, then an almost identical argument will show that $a$
is in the orbital $(bf_0, df_0)$ of $f_1^{f_0}$, the orbital $(bf_0^2,
df_0^2)$ of  $f_1^{f_0^2}$, and also in an orbital of $f_0f_1^{-1}$.  Thus by
Lemma \ref{commutedisj}, these three orbitals must all be equal, which
contradicts the fact that  $(bf_0, df_0) \not= (bf_0^2, df_0^2)$ .   $\square$

\bl \label{nested}   Let $f_0,f_1 \in \ploi$ \ such that $[f_1^{f_0},f_0f_1^{-1}] = 1$ and \ $[f_0f_1^{-1},f_1^{f_0^2}] = 1$.  If $A$ is an orbital of $f_0$ then either there exists an orbital $B$ of $f_1$ where $B \subseteq A$, or $A$ is disjoint from the support of $f_1$.

\el

{\em{Proof:}}  By Lemma \ref{nests},  if  $A \cap B \not= \emptyset$, then either  $A \subseteq B$ or $B \subseteq A$.  Suppose $A \subsetneqq B$.  Let $B = (b,d)$ and let $A = (a,c)$.   Suppose $b$ is not a fixed point of
 $f_0$. Then it must be in some orbital $(r,s)$ of $f_0$, where $r<b<s<a$.
 This contradicts Lemma \ref{nests}. It is similarly shown that $d$ is a fixed point of $f_0$.  So then $(bf_0, df_0) = (b,d) = B$ is an orbital of $f_0^{-1}f_1f_0 = f_2$.  By Remark \ref{infinite}, $[f_2^{f_1},f_1f_2^{-1}] = 1$ and \ $[f_1f_2^{-1},f_2^{f_1^2}] = 1$.  Then by Lemma \ref{eqcommute}, since $B$ is an orbital for both $f_1$ and $f_2$, we have that $f_1$ and $f_2$ must commute in $B$.

So $f_1f_2|_B = f_2f_1|_B$.  Then $f_2|_B = f_1^{-1}f_2f_1|_B = f_3|_B =
f_0^{-1}f_2f_0|_B$.  So then $f_2f_0|_B = f_0f_2|_B$.  So, by Lemma
\ref{commutedisj},  either the orbitals $A$ of $f_0$ and $B$ of $f_2$ must
either be disjoint or equal.  But we assumed that $A \subsetneqq B$, which
gives a contradiction.  $\square$

\subsection{Necessary Conditions}

\bl \label{nicebumpup} Let $f_0,f_1 \in \ploi$ \ such that $[f_1^{f_0},f_0f_1^{-1}] = 1$ and \ $[f_0f_1^{-1},f_1^{f_0^2}] = 1$.   Let $A = (a,c)$ be an orbital of $f_0$ that is an up-bump.  Suppose that $(a,c)$ is not an orbital of $f_1$.  Let $(b_1, d_1)$, $(b_2, d_2) , \dots , (b_n, d_n)$ be all of the orbitals of $f_1$, in increasing order, where $(b_i, d_i) \cap (a,c) \not= \emptyset$.  Then all of the following are true.\\
i. $a < b_1$.\\
ii. There is a point $p<c$ such that $f_0|_{[p,c]} = f_1|_{[p,c]}$.\\
iii.  $d_n = c$.\\
iv. If $p$ is the minimal point where $f_0|_{[p,c]} = f_1|_{[p,c]}$, then $b_nf_0 \geq p$.\\
v. $b_1f_0 > b_n$.

\el

{\em{Proof:}}  Recall from Lemma \ref{nested} that if $B$ is an orbital of $f_1$ and $(a,c)$ is an orbital of $f_0$, then either $B \subseteq (a,c)$ or $B \cap (a,c) = \emptyset$.  So then $a \leq b_1$ and $d_n \leq c$.

{\em{i.}} Suppose that $a = b_1$.  Since $(a,c)$, is not an orbital of $f_1$,
we have $d_1 \not= c$.  Since $f_0$ is an up-bump, $(a,d_1) \subset (a, d_1f_0) \subset (a, d_1f_0^2)$, where $(a,d_1f_0)$ is an orbital for $f_2$ and $ (a, d_1f_0^2)$ is an
orbital for $f_3$. Recall that $f_2 = f_1^{f_0} $,  $f_3 = f_2^{f_0}$, and $[f_2,
f_0f_1^{-1}] = [f_3,
f_0f_1^{-1}] = 1  $.  Since $d_1f_1f_0^{-1} = d_1f_0^{-1} < d_1$, we see that $d_1$ is in some orbital of
$f_0f_1^{-1}$.  So by
Lemma \ref{commutedisj}, both $(a, d_1f_0)$ and $(a, d_1f_0^2)$ must be
orbitals for $f_0f_1^{-1}$.  These are non-disjoint, non-equal intervals, thus
they can not both be orbitals of  $f_0f_1^{-1}$, a contradiction.

{\em{ii.}}  Suppose not.  Then
there is a number $t \in (a,c)$ so that $(t,c)$ is an orbital of
$f_0f_1^{-1}$. Since $c$ is an endpoint of an orbital of $f_0$ and $a < b_1$
(from {\em{i}}), it is
possible to find an $k \in \N$ so that $t < b_1f_0^k < d_1f_0^k \leq c$. \ \
$(b_1f_0^k, d_1f_0^k)$ is an orbital for $f_{k+1}$, which commutes with
$f_0f_1^{-1}$.  So then by Lemma \ref{commutedisj}, either the orbitals
$(t,c)$ and $(b_1f_0^k, d_1f_0^k)$ must be disjoint or equal, a contradiction.

{\em{iii.}} This follows immediately from {\em{ii}}.

{\em{iv.}}  From {\em{iii}}, $(b_nf_0, c)$ is an orbital for
$f_2$.  Since
$f_0|_{[p,c]} = f_1|_{[p,c]}$, we know $[p,c]$ is disjoint from the orbitals of
$f_0f_1^{-1}$.  Also, since $p$ is minimal, there is a $w \in [a,p)$ so that
$(w,p)$ is an orbital of $f_0f_1^{-1}$.  By Lemma \ref{commutedisj}, it must
be the case that  $(b_nf_0, c)$ is disjoint from the orbitals of
$f_0f_1^{-1}$, giving $b_nf_0 \geq p$.

{\em{v.}}  Suppose that  $ b_1 f_0 \leq  b_n$.  Recall that the first orbital
of $f_2$ inside $(a,c)$ is $(b_1f_0, d_1f_0)$.  Since $f_1f_kf_1^{-1} = f_{k+1}$
for all $k \geq 2$, the first orbital of $f_{k+2}$ contained in $(a,c)$ is
$(b_1f_0f_1^{k}, d_1f_0f_1^{k})$.   Since $f_0f_kf_0^{-1} = f_{k+1}$
for all $k \geq 2$, the first orbital of $f_{k+2}$ contained in $(a,c)$ is
$(b_1f_0f_0^{k}, d_1f_0f_0^{k})$.  From the two previous statements, we see that
for all $k \geq 2$, we have $b_1f_0f_0^{k} = b_1f_0f_1^{k}$.  Since $b_1 \in
(a,c)$, an orbital of $f_0$, then there is a $k \in \N$ where $b_1f_0f_0^k >
b_n$.  Since $b_1f_0 \leq b_n$, and, by {\em{iii}}, $(b_n, c)$ is an orbital
of $f_1$, then for all $k \in \N$, we have $b_1f_0f_1^k \leq b_n$, which the
contradiction that for appropriately large $n$,  $ b_n \geq  b_1f_0f_1^k  =
b_1f_0f_0^k > b_n$ .   $\square$

\vspace{12pt}
Similarly, we get the analogous result for down bumps below.

\bl \label{nicebumpdown}
 Let $f_0,f_1 \in \ploi$ \ such that $[f_1^{f_0},f_0f_1^{-1}] = 1$ and \
 $[f_0f_1^{-1},f_1^{f_0^2}] = 1$.   Let $A = (a,c)$ be an orbital of $f_0$
 that is a down-bump.  Assume that $(a,c)$ is not an orbital of $f_1$.  Let
 $(b_1, d_1)$, $(b_2, d_2), \dots , (b_n, d_n)$ be all of the orbitals of
 $f_1$, in increasing order, where $(b_i, d_i) \cap (a,c) \not= \emptyset $.
 Then all of the following are true.  
i. $d_n < c$.  
ii. There is a point $\rho > a$ such that $f_0|_{[a,\rho]} =
f_1|_{[a,\rho]}$.  
iii.  $b_1 = a$.  
iv. If $\rho$ is the maximal point where $f_0|_{[a,\rho]} = f_1|_{[a,\rho]}$,
then $d_1f_0 \leq \rho$.  
v. $d_nf_0 < d_1$.
\el

\bl \label{rel2} Let $f_0, f_1$ be elements in $\ploi$.   Let $(a,c)$ be an orbital of $f_0$.  Suppose that $(a,c)$ is not an orbital of $f_1$.  Assume that whenever $B$ is an orbital of $f_1$, either $B \subseteq (a,c)$ or $B \cap(a,c) = \emptyset$.   Let $(b_1, d_1), (b_2, d_2), \dots , (b_n, d_n)$ be all of the orbitals of $f_1$, in increasing order, that are properly contained in $(a,c)$.  If $f_0$ and $f_1$ either satisty i - v in Lemma \ref{nicebumpup} or i - v in Lemma \ref{nicebumpdown}, then $[f_0f_1^{-1}, f_1^{f_0^2}]|_{[a,c]} = 1$.
\el

{\em{Proof:}}  Suppose that $f_0$ is an up-bump.  So then $f_0$ and $f_1$ satisty {\em{i - v}} in Lemma \ref{nicebumpup}.   The first orbital of $f_1^{f_0^2}$ is $(b_1f_0^2, d_1f_0^2)$.
By {\em{v}}, $b_1f_0 > b_n$ so then  $b_1f_0f_0 > b_nf_0 \geq p$ (by
{\em{iv}}).  So all of the orbitals of $f_1^{f_0^2}$ that are contained in
$(a,c)$ are in $(p,c)$.  By {\em{ii}}, $f_0|_{[p,c]} = f_1|_{[p,c]}$, so then
$f_0f_1^{-1}|_{[p,c]} = 1$.  Then the orbitals of  $f_1^{f_0^2}|_{[a,c]}$ and
$f_0f_1^{-1}|_{[a,c]}$ are disjoint, thus the functions commute.

A similar argument shows that  $[f_0f_1^{-1}, f_1^{f_0^2}]|_{[a,c]} = 1$ when
$(a,c)$ is a down-bump of $f_0$.   $\square$

\vspace{8pt}

The lemma below characterizes the ``nice'' standard isomorphic copies of $F$,
in that there are only 6 properties to check, all of which are completely
determined by the function values of some key points of the domain.  It is not
a necessary condition that $b_1f_0 \geq p$ (for up-bumps) or  $d_1f_0 \leq
\rho$ (for down-bumps) in order for the functions $f_0$ and $f_1$ to satisfy
the standard $F$ relations.

\bl \label{niceiff}   Let $f_0, f_1$ be non-commuting functions in $\ploi$
where for every orbital $A$ of $f_0$ and $B$ of $f_1$, either $B \subseteq A$ or $A \cap B = \emptyset$.  Assume that for
every orbital $(a,c)$ that is an up-bump of $f_0$, if $(b_1, d_1), (b_2, d_2), \dots , (b_n, d_n)$ are all of the orbitals of $f_1$, in increasing order, that are properly contained in $(a,c)$, then i - v in Lemma \ref{nicebumpup} are
satisfied as well as the additional condition $b_1f_0 \geq p$.  Also assume
that for every orbital  $(a,c)$ that is a down-bump of $f_0$, if $(b_1, d_1),
(b_2, d_2), \dots , (b_n, d_n)$ are all of the orbitals of $f_1$, in increasing
order, that are properly contained in $(a,c)$, then i - v in Lemma \ref{nicebumpdown} are satisfied as well as the additional condition $d_nf_0 \leq \rho$.  Then $f_0$ and $f_1$ generate a standard isomorphic copy of $F$.

\el

{\em{Proof:}}  By Lemma \ref{npnaq}, since $f_0$ and $f_1$ do not commute, it only needs to be shown that $[f_0f_1^{-1}, f_1^{f_0}] = 1$ and $[f_0f_1^{-1}, f_1^{f_0^2}] = 1$. 

Since for every orbital $A$ of $f_0$ and $B$ of $f_1$, either $B \subseteq A$ or $A \cap B = \emptyset$, then  $[f_0f_1^{-1}, f_1^{f_0}]|_{[0,1] - \supp(f_0)} =  [f_1^{-1}, f_1]|_{[0,1] - \supp(f_0)}  = 1$ and $[f_0f_1^{-1}, f_1^{f_0^2}]|_{[0,1] - \supp(f_0)}  =  [f_1^{-1}, f_1]|_{[0,1] - \supp(f_0)}  = 1$.  Thus the relations in (2) only need to be checked in the orbitals of $f_0$.

Let $(a,c)$ be an orbital of $f_0$.  If $(a,c)$ is an up-bump, then $p \leq b_1f_0
< b_1f_0^2 < c$.  Recall that when $f_0$ is an up-bump in $(a,c)$, the part of $\supp(f_1^{f_0})$ contained in the
interval $(a,c)$ is equal to  $(b_1f_0, c)$ and the part of
$\supp(f_1^{f_0^2}) $ that is contained in the interval $(a,c)$ is equal to $
(b_1f_0^2, c)$.  This gives $ \supp(f_1^{f_0^2})|_{(a,c)}  \subseteq  \supp(f_1^{f_0})|_{(a,c)} \subseteq (p,c)$ and
$\supp(f_0f_1^{-1}) \cap (p,c) = \emptyset$.  So we have $[f_1^{f_0}, f_0f_1^{-1} ]|_{[a,c]} =
1 $ and $[f_1^{f_0^2}, f_0f_1^{-1} ]|_{[a,c]} = 1$.

Similarly, when $(a,c)$ is a down-bump of $f_0$, then $a < d_nf_0^2 < d_nf_0
\leq \rho$.  So then we have $ \supp(f_1^{f_0^2})|_{(a,c)}  \subseteq
\supp(f_1^{f_0})|_{(a,c)} \subseteq (a,\rho)$ and $\supp(f_0f_1^{-1}) \cap
(a,\rho) = \emptyset$.  So again we have $[f_1^{f_0^2}, f_0f_1^{-1} ]|_{[a,c]} =
1 $ and $[f_1^{f_0}, f_0f_1^{-1} ]|_{[a,c]} = 1$.  $\square$

\subsection{Necessary and Sufficient Conditions}

Previously in this paper, we have shown that there are several conditions that $f_0$ and $f_1$ must satisfy in order to generate a standard isomorphic copy of $F$.  It has also been shown in Lemma \ref{rel2} that these conditions are sufficient to guarentee that the relation $[f_0f_1^{-1}, f_1^{f_0^2} ] = 1$ holds.  The lemmas in this section will focus on the relation $[f_0f_1^{-1}, f_1^{f_0} ] = 1$.

\bl \label{whereequal}  Assume $f_0$ and $f_1$ generate a standard isomorphic
copy of $F$ in $\ploi$.

 Suppose that $(a,c)$ is an up-bump of $f_0$, $(b_n,c)$ is an orbital of $f_1$, and $p$ is the minimal number so that $f_0|_{[p,c]} = f_1|_{[p,c]}$.   If there is a number $q \in (a, p) $ so that $qf_0 = qf_1$, then $q \in (b_n, p)$.

 Alternately, suppose that $(a,c)$ is a down-bump of $f_0$, $(a,d_1)$ is an orbital of $f_1$, and $\rho$ is the minimal number so that $f_0|_{[a,\rho]} = f_1|_{[a,\rho]}$.   If there is a number $q \in (\rho, c) $ so that $qf_0 = qf_1$, then $q \in (\rho, d_1)$.
\el

{\em{Proof:}}  In the case that $(a,c)$ is an up-bump of $f_0$, assume that  $q \notin (b_n, p)$.  Since  $q$ is in an orbital of
$f_0$, $qf_1 = qf_0 \not= q$.  So then $q$ must also be in some orbital,
$(b_k, d_k)$, of $f_1$.  So then $qf_0 = qf_1 \in (b_k, d_k)$, where $d_k \leq
b_n$.  In particular,  $ qf_0 = qf_1 < d_k \leq b_n$.  Combining the previous
inequality with $b_kf_0 < qf_0$ and Lemma \ref{nicebumpup} part {\em{v}}, we
get $qf_0 = qf_1 <  d_k \leq b_n < b_1f_0 \leq b_kf_0 < qf_0 = qf_1 $, a
contradiction. 

The down-bump case is similar.  $\square$

\bl \label{firstr}

Assume $f_0$ and $f_1$ generate a standard isomorphic
copy of $F$ in $\ploi$.  Suppose that $(a,c)$ is an orbital of $f_0$ that
contains the orbitals $(b_1, d_1)$, ... $(b_n, d_n)$ of $f_1$.

Suppose that $(a,c)$ is an up-bump of $f_0$.  If $r$ is the minimal number in
the interval $(a,c)$ so that $rf_0 = rf_1$, then $b_1f_0 \geq r$.

Alternately, suppose that $(a,c)$ is a down-bump of $f_0$.  If $r$ is the maximal number in
the interval $(a,c)$ so that $rf_0 = rf_1$, then $d_nf_0 \leq r$.

\el

{\em{Proof:}}  In the case that $(a,c)$ is an up-bump of $f_0$, assume that
$r$ is the minimal number where $rf_0 = rf_1$.  From Lemma \ref{whereequal},
we have $r > b_n$.  Since $r$ is minimal, then
$(a, r)$ must be an orbital of $f_0f_1^{-1}$.

The interval $(b_1f_0, d_1f_0)$ is the first orbital of $f_1^{f_0}$ contined
in $(a,c)$.   Since the orbitals of $f_0f_1^{-1}$ and $f_1^{f_0}$ must commute, then their orbitals must be disjoint or equal (Lemma \ref{commutedisj}). 

As $(a,r)$ can not equal the interval $(b_1f_0, d_1f_0)$, it follows
that they must be disjoint, and $b_1f_0 \geq r$.

The case for down-bumps is similar.  $\square$

\vspace{8pt}

Define a {\em{nice orbital}} of $f_0$ to be an orbital of $f_0$ where $f_0$
and $f_1$ satisfy the conditions of Lemma \ref{niceiff}.  For convenaince, we
may refer to {\em{nice up-bumps}} and {\em{nice down-bumps}} of $f_0$. 

Lemma \ref{niceiff} gives that any two functions $f_0$ and $f_1$ where every
orbital of $f_0$ is either a nice orbital (or is an interval where $f_0$ and
$f_1$ commute) will generate a standard isomorphic copy
of $F$.  We will call such a pair a {\em{nice generating pair}}.

We will reuse the notation in Lemmas \ref{nicebumpup},
\ref{whereequal}, and \ref{firstr} to categorize the orbitals of $f_1$ into
three disjoint sets.

\be
\item The {\em{Outside Orbitals,}} which are orbitals contained in the
interval $(rf_0^{-1}, pf_0^{-1})$ when $(a,c)$ is an up-bump, and orbitals
contained in the interval $(\rho f_0 , r f_0)$ when $(a,c)$ is a down-bump.

\item The {\em{Inside Orbitals,}} which are orbitals contained in the
interval $(pf_0^{-1} , b_n)$ when $(a,c)$ is an up-bump, and orbitals
contained in the interval $(d_1, \rho f_0)$ when $(a,c)$ is a down-bump.  

\item The {\em{Main Orbital,}} which is the orbital $(b_n,c)$ when $(a,c)$ is an
    up-bump and $(a, d_1)$ when $(a,c)$ is a down-bump.

\ee

Recall that Lemma \ref{firstr} guarantees that for an up-bump, $(b_1, pf_0^{-1})$ is contained
in the interval $(rf_0^{-1}, pf_0^{-1})$.  Also, for a down-bump, the same
lemma gives that $(\rho f_0 , d_n )$ is contained in  $(\rho f_0 , r f_0)$.  

 If $f_1$ does not have any Outside Orbitals in $(a,c)$, then $(a,c)$ would be
 a nice orbital of $f_0$, and all of the relations of $F$ would be satistied
 on this interval.

\bt \label{mainthm}
Let $f_0$ and $f_1$ be non-commuting functions in $\ploi$ that satisty i - v in 
Lemmas \ref{nicebumpup} (for up-bumps) and  \ref{nicebumpdown} (for down bumps).  Suppose also that every
orbital of $f_1$ that is contained in an orbital of $f_0$ is either a Main
Orbital, an Inside Orbital, or an Outside Orbital.   The functions $f_0$ and $f_1$ generate a
standard isomorphic copy of $F$ if and only if every orbital of $f_1^{f_0}$
that corresponds to an Outside Orbital of $f_1$ commutes with $f_0f_1^{-1}$.
\et

{\em{Proof:}}

$(\Rightarrow)$
We have shown in several previous lemmas that if $f_0$ and $f_1$ generate a
standard isomorphic copy of $F$, then all of the listed properties must be true.

$(\Leftarrow)$
Since $f_0$ and $f_1$ are non-commuting, if they satisfy the relations of $F$
then they generate a standard isomorphic copy of $F$.

As all orbitals of $f_1$ satisfy {\em{i - v}} of \ref{nicebumpup} or
\ref{nicebumpdown}, we have shown that these properties alone are enough to
guarantee the second relation, $[f_0f_1^{-1},
 f_1^{f_0^2}]=1$, holds.

Every orbital of $f_1^{f_0}$ that corresponds to a Main Orbital or an Inside
Orbital of $f_1$ is disjoint from the support of $f_0f_1^{-1}$, thus they
commute.  In addition, since every orbital of $f_1$ that does not commute with
$f_0$ is an Outside Orbital and every orbital of $f_1^{f_0}$ that corresponds to an
Outside Orbital of $f_1$ also commutes with $f_1f_1^{-1}$, then we get the
first relation $[f_0f_1^{-1},
 f_1^{f_0}]=1$.  $\square$

\subsection{Building a Standard Generating Pair}

We will end this section with a description of exactly how to construct any
pair of functions, $f_0$ and $f_1$, that generate a standard isomorphic copy of $F$.

To make an arbitrary pair of functions that generate a standard isomorphic
copy of $F$, we will begin with two functions $f_0$ and $f_1$ where every
orbital of $f_0$ is either a nice orbital or is an interval where $f_0$ and
$f_1$ commute.  We will then slightly modify the functions in a nice orbital
of $f_0$ to construct two new
functions that contain Outside Orbitals of $f_1$ and also generate a standard
isomorphic copy of $F$.

First, we will explicitely state how to build any nice generating pair $f_0$
and $f_1$, where $f_0$ has exactly one nice up-bump.  The process is almost
identical for constructing pairs with exactly one nice down-bump, as well as
constructing pairs with more than one nice orbital. 

\br \label{niceorbital}
Let $f_0$ be any function in $\ploi$ that is not the identity function.  In
any up-bump $(a,c)$ of $f_0$, pick any point $p \in (a,c)$.  Let $f_1$ be any
function in $\ploi$ where $f_0 = f_1$ on the interval $[p,c]$, where $f_1$
is the identity in the interval $(a, pf_0^{-1})$, and where $f_1$ and $f_0$
commute outside of $(a,c)$.  The funcions $f_0$ and
$f_1$ generate a standard isomorphic copy of $F$.
\er

\bt \label{construct1}
Assume that $(a,c)$ is a nice up-bump of $f_0$.  Let $p$ be as in
Lemma \ref{niceiff}.  Let $s \in (p, c)$.  If $h$ be any function in $\ploi$
where $\supp(h) \subseteq (pf_0^{-1}, sf_0^{-1})$ and $t,k$ are any
integers, then the pair $f_0$ and $f_1$ are a standard generating pair of $F$
if and only if $f_0$ and $g_1 = h^tf_0^{-1}h^kf_0f_1$ are a standard
generating pair of $F$.  
\et

{\em{Proof:}}  As the functions $f_1$, and $g_1$ are identical outside of the
interval $(a,c)$, then the pair $f_0$ and $f_1$ satisfy the relations of $F$
outside the interval $(a,c)$ if and only if the relations are also satistied
outside of $(a,c)$ by the pair $f_0$ and $g_1$.  As $f_0$ and $f_1$ satisfy
the conditions of Lemma \ref{niceiff}, then if $f_0$ and $g_1$ are
a standard generating pair, then so are $f_0$ and $f_1$.

Suppose that $f_0$ and $f_1$ are a standard generating pair for $F$.  By
construction, $f_0$ and $g_1$ satisfy the conditions in Lemma \ref{rel2}, thus
only the first relation must be checked for these functions.  Using the fact
that $f_0f_1^{-1}$ is the identity on the interval $(p,c)$, it is a
straightforward to check that whenever $y \in (p,s)$, we have $yf_0g_1^{-1}g_1^{f_0} =
yg_1^{f_0}f_0g_1^{-1} = yf_0^{-1}h^{t-k}f_0$.  Combine this computation with
the facts that $g_1^{f_0}$ is the
identity in the interval $(a,p)$ and $f_0g_1^{-1}$ is the identity in the
interval $(s,c)$ to get that the relation $[g_1^{f_0}, f_0g_1^{-1}]=1$ holds.  $\square$

\vspace{8pt}

\bt \label{construct2}
If $f_0$ and $g_1$ are a standard generating pair of $F$, then $f_0$ and $g_1$
can be inductively constructed from a nice generating pair $f_0$, $f_1$ in a
finite number of steps as described in Theorem \ref{construct1} (or the
theorem's analogous result for down-bumps).   
\et

{\em{Proof:}}  Assuming that $f_0$ and $g_1$ are not a nice generating
pair of $F$, there must be an orbital $(a,c)$ of $f_0$ that properly contains an
Outside Orbital of $g_1$. Without loss of
generality, assume $(a,c)$ is an up-bump.  Let $p$ be the minimal point where
$f_0|_{[p,c]} = g_1|_{[p,c]}$.  Let $\alpha$ be the minimal point in $(a,c)$
where $\alpha f_0 = \alpha g_1$.  The first orbital $(b_1f_0, d_1f_0)$ of $g_1^{f_0}$ in $(a,c)$ must either be an orbital
of $f_0g_1^{-1}$ or disjoint from the support of $f_0g_1^{-1}$.  In either
case, $(b_1f_0)f_0 = (b_1f_0)g_1$, so then $\alpha \leq b_1f_0 < p $.    From
Lemma \ref{whereequal}, $b_n < \alpha$. 

Define $f_1$ to be the function in $\ploi$ that is
equal to $g_1$ outside of the interval $(a, c)$, is the identity
function in the interval $(a, pf_0^{-1})$, is equal to $g_1$ in the
interval $( pf_0^{-1} , \alpha  )$, and is equal to $f_0$ in the
interval $(\alpha , c)$.  As all the conditions of Lemma \ref{niceiff} are
satistied for $f_0$ and $f_1$, the pair $f_0$ and $f_1$ are a standard
generating pair for $F$, and has one less non-nice orbital than $f_0$ and
$g_1$.

Let $(b_1, d_1), ..., (b_m, d_m)$ be the orbitals of $g_1$ that are not
orbitals of $f_1$.  Let $(q_1, w_1), ... , (q_z, w_z)$ be the orbitals of
$f_0g_1^{-1}$ that are contained in $(\alpha, c)$.  As there are no non-equal,
overlapping intervals in the set $S = \{ (b_if_0, d_if_0) , (q_j, w_j)  \}$, we will
always look at the first orbital in the set.  This will give three cases for
each step in the inductive process.  Let $g_{(1,0)} = f_1$.  The following
three cases will construct $g_{(1,1)}$.

For the first case, assume $(b_1f_0, d_1f_0)$ is the unique first interval in
the set $S$.  Then let $h_1$ be the function in $\ploi$ where $\supp(h_1)
\subseteq (b_1, d_1)$ and  $h_1|_{(b_1, d_1)} = g_1|_{(b_1, d_1)}$.  Let $g_{(1,1)} =
h_1^1f_0^{-1}h_1^0f_0g_{(1,0)} = h_1f_1$, which equals $g_1$ on $(b_1, d_1)$ and
$(pf_0^{-1}, d_1f_0)$.  Now remove $(b_1f_0, d_1f_0)$ from the set $S$ and relabel
$(b_if_0, d_if_0)$ in $S$ to $(b_{i-1}f_0, d_{i-1}f_0)$. 

For the second case, assume that $(q_1, w_1)$ is the unique first interval in
the set $S$.  Then let $\psi$ be the function in $\ploi$ where $\supp(\psi)
\subseteq (q_1, w_1)$ and $\psi|_{(q_1, w_1)} = f_0g_1^{-1}|_{(q_1, w_1)}$.
Let $h_1 = f_0 \psi^{-1} f_0^{-1}$.  Then if we let $g_{(1,1)} =
h_1^0f_0^{-1}h_1^1f_0g_{(1,0)}$, we have $g_{(1,1)}$ equals $g_1$ on $(q_1f_0^{-1},
w_1f_0^{-1})$ and $(pf_0^{-1}, w_1)$.  Now remove $(q_1, w_1)$ from the set
$S$ and relabel $(q_j, w_j)$ in $S$ to $(q_{j-1}, w_{j-1})$.

For the third case, assume $(b_1f_0, d_1f_0) = (q_1, w_1)$.  Then by Theorem \ref{commonroot}, $g_1^{f_0}|_{(b_1f_0, d_1f_0)}$ and
$f_0g_1^{-1}|_{(b_1f_0, d_1f_0)}$ are powers of a common root function.  Let
$\phi$ be the function in $\ploi$ whose support is contained in $(b_1f_0,
d_1f_0)$ and $t$ and $-k$ be the integers so that on the interval  $(b_1f_0,
d_1f_0)$, we have $\phi^t = g_1f_0$ and $\phi^{-k} = f_0g_1^{-1}$.  In
this case, let $h_1 = \phi^{f_0^{-1}}$.  Then the function $g_{(1,1)} =
h_1^tf_0^{-1}h_1^kf_0g_{(1,0)}$ is equal to $f_1$ on the intervals $(b_1, d_1)$ and
$(pf_0^{-1}, d_1f_0)$.  Now remove both $(b_1f_0, d_1f_0)$ and  $(q_1, w_1)$
from $S$ and relabel $(b_if_0, d_if_0)$ to $(b_{i-1}f_0, d_{i-1}f_0)$  and $(q_j, w_j)$ to $(q_{j-1}, w_{j-1})$.

In each of the three cases, $S$ becomes smaller by one or two intervals.  For
the $\beta^{th}$ step in the inductive process, find the function
$h_{\beta}$ and integers $r,k$ as in steps
1 -- 3 above, and let $g_{(1, \beta)} = h_1^tf_0^{-1}h_1^kf_0g_{(1,\beta -
  1)}$.  After a finite amount of steps, $S$ will become empty and the last
function created, $g_{(1, \sigma)}$, will be equal to $g_1$.   $\square$

\section{Partial Converse of the Ubiquity Result}

The following theorem is the Ubiquity result, which is the main theorem in \cite{BrinU}.

\vspace{10pt}

\noindent
{\bf{Theorem 3.1}}  {\em{If $H$ is a subgroup of $\ploi$ and there exists an orbital $W = (a,b)$ of $H$ and a function $f$ in $H$ so that either $f|_W$ moves points near $a$ but not $b$ or $f|_W$ moves points near $b$ but not $a$, then there is a subgroup $G$ of $H$ so that $G$ is isomorphic to $F$.  }}

\vspace{10pt}

A partial converse to the Ubiquity result has previously been proven by Guba in an unpublished paper, whose proof used much more machinery than what was used in this thesis.  A proof with significantly less machinery is given below.

\vspace{10pt}

\noindent
{\bf{Theorem 3.2}}  {\em{If $F \cong H \leq \ploi$, then there exists an orbital $W = (a, c)$ of $H$ and a function $h$ in $H$ where  either $h|_W$ moves points near $a$ but not $c$ or $h|_W$ moves points near $c$ but not $a$.}}

\vspace{10pt}

{\em{Proof:}}  Since  $F \cong G$, there are functions $f_0$ and $f_1$ in $H$
that generate a standard isomorphic copy of $F$.  Then by Theorem
\ref{geniff}, there exists an orbital $(a,c)$ of $f_0$ that properly contains
an orbital of $f_1$.  From Theorem \ref{geniff} and Lemma \ref{unionSupports},
we see that $(a,c)$ is an orbital $W$ of $H$.  By Lemma \ref{geniff}, if
$(a,c)$ is an up-bump of $f_0$ then $f_1$ moves points near $c$ but not $a$,
and if $(a,c)$ is a down-bump of $f_0$ then $f_1$ moves points near $a$ but
not $c$.  $\square$

\vspace{8pt}

The full converse of the Ubiquity result is false.   That is, if $F \cong G \leq H \leq \ploi$, then we cannot reach a conclusion about $H$ as given in Theorem 3.2.

For a counterexample, let $x_0$ and $x_1$ be the standard generators of $F$
defined in the introduction.  Recall that for all $k > 0$,  $x_{k+1} =
x_1^{x_0^k}$.   Consider the group $H = \langle x_0, f_0 , f_1 \rangle $,
where $f_0 = x_1^{2}x_2^{-1}x_1^{-1}$ and $f_1 = x_1
x_2^2x_3^{-1}x_2^{-1}x_1^{-1} $.

From Lemma \ref{niceiff}, $\langle f_0, f_1 \rangle \cong F$.  From our
knowledge about the support of a subgroup (Lemma \ref{unionSupports}), we have
$(0, 1) = \supp(x_0) \subseteq \supp(H) \subseteq (0, 1)$, which makes $(0,1)$
the only orbital of $H$.  Let $h$ be any element in $H$.  Since both $f_0$ and
$f_1$ are the identity near $0$, there is an $\varepsilon \in (0, 1)$ and a
unique $n \in \Z$ so that $h|_{(0, \varepsilon)} = x_0^n|_{(0, \varepsilon)}$.
Similarly, there is a $\delta \in (0,1)$ and a unique $m \in \Z$ so that $h|_{(\delta, 1)} = x_0^m|_{(\delta, 1)}$.  Since the support of $x_0$ is $(0,1)$, it must be the case that $n = m$.  So then $h|_{(\delta, 1)} $ is the identity if and only if $m = n = 0$ if and only if $h|_{(0, \varepsilon)}$ is the identity.  Thus $H$ has no element that moves points near $0$ but not near $1$, nor vice versa.

\def\cprime{$'$}
\providecommand{\bysame}{\leavevmode\hbox to3em{\hrulefill}\thinspace}
\providecommand{\MR}{\relax\ifhmode\unskip\space\fi MR }
\providecommand{\MRhref}[2]{%
  \href{http://www.ams.org/mathscinet-getitem?mr=#1}{#2}
}
\providecommand{\href}[2]{#2}

\end{document}